\documentclass[journal]{IEEEtran}
\usepackage{amssymb,latexsym,amsmath}
\usepackage{epsfig}
\def\Ebox#1#2{%
\medskip
\begin{center}
  \strut\epsfxsize=#1 \hsize\epsfbox{#2}
\end{center}
\smallbreak}
\usepackage{caption}

\def\mindex#1{\index{#1}}

\def\sq{\hbox{\rlap{$\sqcap$}$\sqcup$}}
\def\qed{\ifmmode\sq\else{\unskip\nobreak\hfil
\penalty50\hskip1em\null\nobreak\hfil\sq
\parfillskip=0pt\finalhyphendemerits=0\endgraf}\fi\medskip}

\long\def\defbox#1{\framebox[.9\hsize][c]{\parbox{.85\hsize}{%
\parindent=0pt
\baselineskip=12pt plus .1pt      
\parskip=6pt plus 1.5pt minus 1pt 
 #1}}}

\long\def\beginbox#1\endbox{\subsection*{}%
\hbox{\hspace{.05\hsize}\defbox{\medskip#1\bigskip}}%
\subsection*{}}

\def\endbox{}

\newsavebox{\junk}
\savebox{\junk}[1.6mm]{\hbox{$|\!|\!|$}}

\def\limsup{\mathop{\rm lim\ sup}}

\def\state{{\sf X}}

\def\bx{{{\cal B}(\state)}}

\newcommand{\field}[1]{\mathbb{#1}}

\def\Re{\field{R}}

\def\ind{\field{I}}

\def\bfmath#1{{\mathchoice{\mbox{\boldmath$#1$}}%
{\mbox{\boldmath$#1$}}%
{\mbox{\boldmath$\scriptstyle#1$}}%
{\mbox{\boldmath$\scriptscriptstyle#1$}}}}

\def\bfmY{\bfmath{Y}}

\def\bfmhhaY{\bfmath{\hhaY}} 
\def\bfmhhaY{\hbox to 0pt{$\widehat{\bfmY}$\hss}\widehat{\phantom{\raise 1.25pt\hbox{$\bfmY$}}}}

\def\til={{\widetilde =}}

\def\clB{{\cal B}}

\def\clF{{\cal F}}

\def\clM{{\cal M}}

 \def\FRAC#1#2#3{\genfrac{}{}{}{#1}{#2}{#3}}

\def\ddtp{{\mathchoice{\FRAC{1}{d^{\hbox to 2pt{\rm\tiny +\hss}}}{dt}}%
{\FRAC{1}{d^{\hbox to 2pt{\rm\tiny +\hss}}}{dt}}%
{\FRAC{3}{d^{\hbox to 2pt{\rm\tiny +\hss}}}{dt}}%
{\FRAC{3}{d^{\hbox to 2pt{\rm\tiny +\hss}}}{dt}}}}

\def\half{{\mathchoice{\FRAC{1}{1}{2}}%
{\FRAC{1}{1}{2}}%
{\FRAC{3}{1}{2}}%
{\FRAC{3}{1}{2}}}}

\def\eqdef{\mathbin{:=}}

\def\Prob{{\sf P}}

\def\Expect{{\sf E}}

\def\average#1,#2,{{1\over #2} \sum_{#1}^{#2}}

\def\eye(#1){{\bf(#1)}\quad}

\def\varble{\,\cdot\,}

\newtheorem{theorem}{Theorem}[section]

\def\Theorem#1{Theorem~\ref{#1}}

\def\Section#1{Section~\ref{#1}}

\def\eq#1/{(\ref{e:#1})}

\newcommand{\beqn}[1]{\notes{#1}%
\begin{eqnarray} \elabel{#1}}

\newcommand{\eeqn}{\end{eqnarray} }

\newcommand{\beq}[1]{\notes{#1}%
\begin{equation}\elabel{#1}}

\newcommand{\eeq}{\end{equation}}

\def\bdes{\begin{description}}
\def\edes{\end{description}}

\newcounter{rmnum}
\newenvironment{romannum}{\begin{list}{{\upshape (\roman{rmnum})}}{\usecounter{rmnum}
\setlength{\leftmargin}{14pt}
\setlength{\rightmargin}{8pt}
\setlength{\itemindent}{-1pt}
}}{\end{list}}

\newcounter{anum}

{\end{list}}

\def\ass(#1:#2){(#1\ref{#1:#2})}

\def\ritem#1{
\item[{\sf \ass(\current_model:#1)}]
}

\newenvironment{recall-ass}[1]{%
\begin{description}
\def\current_model{#1}}{
\end{description}
}

\def\Ebox#1#2{%
\begin{center}
 \parbox{#1\hsize}{\epsfxsize=\hsize \epsfbox{#2}}
\end{center}}

\newcommand{\bd}{\begin{description}}
\newcommand{\ed}{\end{description}}
\newcommand{\bt}{\begin{theorem}}
\newcommand{\et}{\end{theorem}}
\newcommand{\ba}{\begin{array}{rcl}}
\newcommand{\ea}{\end{array}}

\def\tn{\zeta}

\def\stp{{\cal T}}

\newlength{\noteWidth}
\long\def\notes#1{\ifinner
           {\tiny #1}
           \else
           \marginpar{\parbox[t]{\noteWidth}{\raggedright\tiny #1}}
       \fi\typeout{#1}}

       \newtheorem{thm}{\bf{Theorem}}[section]
       \newtheorem{cor}{\bf{Corollary}}[section]
       \newtheorem{lem}{\bf{Lemma}}[section]
       
       \newtheorem{prop}{\bf{Proposition}}[section]
       \newtheorem{defn}{\bf{Definition}}[section]
 \newtheorem{remark}{\bf{Remark}}[section]
\begin{document}

\title{Random-Time, State-Dependent Stochastic Drift for Markov Chains and Application to Stochastic Stabilization Over Erasure Channels$^1$}
\author{Serdar Y\"uksel$^2$ and Sean P. Meyn$^3$}
\maketitle
\footnotetext[1]{The material of this paper was presented in part at the 2009 Annual Allerton Conference on Communication, Control, and Computing, in September 2009.}
\footnotetext[2]{Department of Mathematics and Statistics, Queen's University, Kingston, Ontario, Canada, K7L 3N6. Email: yuksel@mast.queensu.ca. Research Supported By the Natural Sciences and Engineering Research Council of Canada (NSERC).}
\footnotetext[3]{Department of Electrical and Computer Engineering, University of Florida, FL 32611 USA. Email: meyn@ece.ufl.edu.}

\maketitle
\begin{abstract}
It is known that state-dependent, multi-step Lyapunov bounds lead to greatly simplified verification theorems for stability for large classes of Markov chain models. This is one component of the ``fluid model'' approach to stability of stochastic networks.   In this paper we extend the general theory to randomized
multi-step Lyapunov theory to obtain criteria
for stability and steady-state performance bounds, such as finite moments.

These results are applied to  a remote stabilization problem, in which a controller receives measurements from
an erasure channel with limited capacity.    Based on the general results in the paper it is shown that stability of the closed loop system is assured provided that the channel capacity is greater than the logarithm of the unstable eigenvalue, plus an additional correction term. The existence of a finite second moment in steady-state is established under additional conditions.
\end{abstract}
\section{Introduction}

Stochastic stability of Markov chains has a rich and complete theory, and forms a foundation for several other general techniques such as dynamic programing and Markov Chain Monte-Carlo (MCMC) \cite{CTCN}.   This paper concerns extensions and application of a class of Lyapunov techniques, known as \textit{state-dependent}
drift
criteria  \cite{MeynStateDrift}.  This technique is the basis of the fluid-model (or ODE) approach to stability in stochastic networks and other general models \cite{bormey00a,faymalmen95a,malmen79,formeymoupri08a,dai95a,daimey95a,CTCN,MCMC}.


In this paper we consider a stability criterion based on a state-dependent \textit{random sampling} of the Markov chain of the following form:
It is assumed that there is a function $V$ on the state space taking positive values, and an increasing sequence of stopping times $\{\stp_i : i \in \mathbb{N}_+\}$, with $\stp_0=0$, such that for each $i$,
\begin{equation}
\Expect[V(x_{\stp_{i+1}}) \mid \clF_{\stp_i}]   \le V(x_{\stp_i})  - \delta(x_{\stp_i})
\label{e:RandomDrift}
\end{equation}
where the function $\delta\colon\state\to \Re$ is positive (bounded away from zero) outside of a ``small set'',
and $\clF_{\stp_i}$ denotes the filtration of ``events up to time $\stp_i$''.
Under suitable conditions on the Markov chain, the drift $\delta $, and the sequence $\{\stp_i\}$, we establish corresponding stability and ergodicity properties of the chain.
The main results extend and unify previous research on stability and  criteria for finite moments obtained in
\cite{malmen79,MeynStateDrift,formeymoupri08a,daimey95a,CTCN}.

Motivation for this research arose from our interest in applications to networked control, and information theory with variable length and variable delay decoding \cite{Burnashev}, \cite{SahaiParts}, \cite{Sahai}, and non-asymptotic information theory \cite{PolyanskiyISIT10}. Specifically, in some
network protocol, team decision and networked control applications there is only intermittent access to sensor information, control action or some common knowledge on which decisions can be made. The timing may be random, depending on availability of communication resources. One example of such conditions is reported in \cite{YukTAC2010}, for establishing stochastic stability of adaptive quantizers for Markov sources where random stopping times are the instances when the encoder can transmit granular information to a controller. We will also consider such an application in detail in the paper. In this context, there has been a significant amount of research on stochastic stabilization of networked control systems under information constraints. For a detailed review see \cite{YukBasTACNoisy}. Stochastic stability of adaptive quantizers have been studied both in the information theory community (see \cite{GoodmanGersho}, \cite{KiefferDunham}) as well as control community (\cite{BrockettLiberzon}, \cite{LiberzonNesic}, \cite{NairEvans}, \cite{YukTAC2010}). \cite{NairEvans} provided a stability result under the assumption that a quantizer is variable-rate for systems driven by noise with unbounded support for its probability measure. \cite{NairEvans} used asymptotic quantization theory to obtain a time-varying scheme, where the quantizer is used at certain intervals at a very high rate, and at other time stages, the quantizer is not used. For such linear systems driven by unbounded noise, \cite{YukTAC2010} established ergodicity, under fixed-rate constraints, through martingale methods. These papers motivated us to develop a more general theory for both random-time drift as well as the consideration of more general noise models on the channels. In a similar line of work, \cite{GurtNair} also considered stability of the state and quantization parameters, \cite{Minero} studied the problem concerning time-varying channels and provided a necessity and sufficiency result for boundedness of second moments, \cite{YouXie} studied the problem of control over an erasure channel in the absence of noise, and \cite{YukITA2011} considered discrete noisy channels with noiseless feedback for systems driven by Gaussian noise using the machinery developed in the current paper. \cite{YukBasTACNoisy} obtained conditions for the existence of an invariant probability measure for noisy channels, considering deterministic, state-dependent drift, based on the criteria developed in \cite{MeynStateDrift}.

We believe our results will provide constructive tools to address related stability issues in a large class of networked control problems.

The contributions of this paper can be summarized as follows

\begin{itemize}
\item Stochastic stability theory for Markov chains based on random-time, state-dependent stochastic drift criteria.   A range of conditions are used to establish conditions for positive recurrence, and existence of finite moments.

\item The results are applied to stochastic stabilization over an erasure network, where a linear system driven by Gaussian noise is controlled. This paper establishes that in such an application of stabilization of an unstable system driven by Gaussian noise over erasure channels, for stochastic stability, Shannon capacity is sufficient (up to an additional correction term). For the existence of finite moments, however, more stringent criteria are needed. Regarding information rate requirements, our construction is tight up to an additional symbol, in comparison with necessary conditions presented in \cite{Minero}.
\end{itemize}

The remainder of the paper is organized as follows.  The implications of the drift criterion \eqref{e:RandomDrift} to various forms of stochastic stability are presented in the next section. The rest of the paper focuses on an application to control over a lossy erasure network with quantized observations.  The paper ends with concluding remarks in Section~\ref{s:conc}.
Proofs of the stochastic stability results and other technical results are contained in the appendix.

\section{Stochastic Stability}
\label{s:ss}
%
\subsection{Preliminaries}

We let ${\bf \phi} =\{ \phi_t, t\geq 0\}$ denote a Markov chain with state space $\state$.
The basic assumptions of \cite{MCSS} are adopted:  It is assumed that $\state$ is a complete separable metric space,  that is locally compact; its  Borel $\sigma$-field is denoted $\clB(\state)$.  The transition probability is denoted by $P$, so that for any $\phi\in\state$,  $A\in\bx$,  the probability of moving in one step from the state $\phi$ to the set $A$  is given by  $  \Prob(\phi_{t+1}\in A \mid \phi_t=\phi) = P(\phi,A)$.   The $n$-step transitions are obtained via composition in the usual way, $  \Prob(\phi_{t+n}\in A \mid \phi_t=\phi) = P^n(\phi,A)$, for any $n\ge1$.   The transition law acts on measurable functions $f\colon\state\to\Re$ and measures $\mu$ on $\bx$ via,
\[
Pf\, (\phi)\eqdef \int_{\state} P(\phi,dy) f(y),\quad \phi\in\state, \]
and
\[ \mu P\, (A) \eqdef \int_{\state} \mu(d\phi)P(\phi,A) ,\quad A\in\bx.
\]
A probability measure $\pi$ on $\bx$ is called invariant if $\pi P= \pi$.  That is,
\[
\int \pi(d\phi) P(\phi,A) = \pi(A),\qquad A\in\bx.
\]


For any initial probability measure $\nu$ on $\bx$ we can construct a stochastic process with transition law $P$, and satisfying $\phi_0\sim \nu$.  We let $\Prob_\nu$ denote the resulting probability measure on sample space, with the usual convention for $\nu=\delta_{\phi}$ when the initial state is $\phi\in\state$. When $\nu=\pi$ then the resulting process is stationary.


There is at most one stationary solution under the following irreducibility assumption.  For a set $A\in\bx$ we denote,
\begin{equation}
\tau_A\eqdef \min(t \ge 1 :  \phi_t \in A)
\label{e:tauA}
\end{equation}
\begin{defn}
Let  $\varphi$ denote  a sigma-finite measure on $\bx$.
\begin{romannum}
\item
The Markov chain is called \textit{$\varphi$-irreducible} if for any $\phi\in\state$,
and  any $B\in\bx$ satisfying $\varphi(B)>0$,     we have
\[
\Prob_{\phi}\{\tau_B<\infty\} >0 \, .
\]

\item
A $\varphi$-irreducible Markov chain is \textit{aperiodic}  if  for any $\phi\in\state$, and any $B\in\bx$ satisfying $\varphi(B)>0$,   there exists $n_0=n_0(\phi,B)$ such that
\[
 P^n(\phi,B)>0 \qquad \hbox{\it for all \ } n\ge n_0.
\]

\item
A $\varphi$-irreducible Markov chain is \textit{Harris recurrent}  if   $\Prob_{\phi}(\tau_B < \infty  ) =1 $ for any $\phi\in\state$, and any $B\in\bx$ satisfying $\varphi(B)>0$.  It is \textit{positive Harris recurrent} if in addition there is an invariant probability measure $\pi$.
\end{romannum}
\end{defn}

Intimately tied to $\varphi$-irreducibility is the existence of  a suitably rich collection of  ``small sets'', which allows Nummelin's splitting technique to be applied leading to verification for Harris recurrence.   A set $A \in\bx$ is small if there is an integer $n_0\ge 1$  and a positive measure $\mu$ satisfying $\mu(\state)>0$ and
$$
P^{n_0}(x,B)  \geq \mu(B), \quad \hbox{\it for all \ }  x \in A, \; \mbox{and} \, B \in \bx .
$$
Small sets are analogous to compact sets in the stability theory for $\varphi$-irreducible Markov chains.
In most applications of  $\varphi$-irreducible Markov chains we find that any compact set is small -- In this case, ${\bf \phi}$ is called a \textit{T-chain} \cite{MCSS}.
%

To relax the $\varphi$-irreducibility assumption we can impose instead the following continuity assumption: A Markov chain is (weak) Feller if  the function $Pf$ is continuous on $\state$, for every continuous and bounded function $f\colon\state\to\Re$.

We next introduce criteria for positive Harris recurrence for $\varphi$-irreducible Markov chains, and criteria for the existence of a steady-state distribution $\pi$ for a Markov chain satisfying the Feller property.

\subsection{Drift criteria for positivity}
\label{driftCriteria}

We now consider specific   formulations of the random-time drift criterion \eqref{e:RandomDrift}.   Throughout the paper the sequence of stopping times $\{\stp_i : i \in \mathbb{N}_+\}$ is assumed to be non-decreasing, with $\stp_0=0$.
In prior work on state-dependent criteria for stability it is assumed that the stopping times take the following form,
\[
\stp_{i+1}=\stp_i + n(\phi(\stp_i)),\qquad i\ge 0
\]
where $n\colon\state\to\mathbb{N}$ is a deterministic function of the state. The results that follow generalize state dependent drift results in \cite{MeynStateDrift} to this random-time setting. We note that a similar approach has been presented recently in the literature in \cite{Fralix} for random-time drift (see Theorem 4), which readily generalizes the state dependent drift results in \cite{MeynStateDrift}. The conditions presented in \cite{Fralix} are more restrictive for the stopping times than what we present here. Furthermore we present discussions for existence of finite moments, as well as extensions for non-irreducible chains.

The proofs of these results are presented in the appendix.

%
%

\Theorem{thm5} is the main general result of the paper, providing a single criterion for positive Harris recurrence, as well as   finite  ``moments'' (the steady-state mean of the function $f$ appearing in the drift condition \eqref{e:thm5delta}).
The drift condition \eqref{e:thm5delta} is a refinement of \eqref{e:RandomDrift}. 

\begin{thm}
\label{thm5}
Suppose that ${\bf \phi}$ is a $\varphi$-irreducible and aperiodic Markov chain.   Suppose moreover that there are functions
 $V \colon \state \to (0,\infty)$,
 $\delta \colon \state \to [1,\infty)$,
 $f\colon \state \to [1,\infty)$,
a small set $C$,   and a constant $b \in \Re$,   such that the following hold:
\begin{equation}
\begin{aligned}
\Expect[V(\phi_{\stp_{z+1}}) \mid \clF_{\stp_z }]  &\leq  V(\phi_{\stp_{z}}) -\delta(\phi_{\stp_z }) + b1_{\{\phi_{\stp_z} \in C\}}
\\
 \Expect \Bigl[\sum_{k=\stp_z}^{\stp_{z+1}-1} f(\phi_k)  \mid \clF_{\stp_z }\Bigr]  &\le \delta(\phi_{\stp_z})\, , \qquad \qquad \qquad \qquad z\ge 0.
\end{aligned}
\label{e:thm5delta}
\end{equation}
Then the following hold:
\begin{romannum}
\item
${\bf \phi}$ is positive Harris recurrent, with unique invariant distribution $\pi$
\item $\pi(f)\eqdef \int f(\phi)\, \pi(d\phi) <\infty$
\item  For any function $g$ that is bounded by $f$, in the sense that $\sup_{\phi} |g(\phi)|/f(\phi)<\infty$, we have convergence of moments in the mean, and the Law of Large Numbers holds:
\[
\begin{aligned}
\lim_{t\to\infty} \Expect_{\phi}[g(\phi_t)] &= \pi(g)
\\
\lim_{N\to\infty} \frac{1}{N} \sum_{t=0}^{N-1} g(\phi_t) &= \pi(g)\qquad a.s.\,, \ \phi\in\state
\end{aligned}
\]
\end{romannum}
\qed
\end{thm}

\begin{remark} We note that, for (ii) in Theorem \ref{thm5}, the condition that  $f \colon \state \to [1,\infty)$,
 $\delta \colon \state \to [1,\infty)$, can be relaxed to $f \colon \state \to [0,\infty)$,
 $\delta \colon \state \to [0,\infty)$, provided that one can establish (i), that is the positive Harris recurrence of the chain first.
\end{remark}

This result has a corresponding, albeit weaker, statement for a Markov chain that is Feller, but not necessarily satisfying the irreducibility assumptions:

\begin{thm}\label{thm25Feller}
Suppose that ${\bf \phi}$ is a Feller Markov chain, not necessarily  $\varphi$-irreducible.
If in addition \eqref{e:thm5delta} holds with $C$ compact,
then there exists at least one invariant probability measure.  Moreover, there exists $c < \infty$ such that, under any invariant probability measure $\pi$,
\begin{equation}
\Expect_{\pi}[f(\phi_t)] = \int_{\mathbb{X}} \pi(d\phi) f(\phi) \leq c.
\label{e:thm5deltaFellerBdd}
\end{equation}
\qed
\end{thm}

We conclude by stating a simple corollary to Theorem~\ref{thm5}, obtained by taking $f(\phi)=1$ for all $\phi \in \state$.
\begin{cor}\label{corol}
Suppose that ${\bf \phi}$ is a $\varphi$-irreducible Markov chain.   Suppose moreover that there is a function $V: \state \to (0,\infty)$, a small set $C$, and a constant $b \in \Re$,   such that the following hold:
\begin{equation}\label{PositiveRecur}
\begin{aligned}
\Expect[V(\phi_{\stp_{z+1}}) &\mid \clF_{\stp_z }] \leq V(\phi_{\stp_z }) - 1 + b1_{\{\phi_{\stp_z} \in C\}}
\\
\sup_{\ z\ge 0} \Expect[\stp_{z+1} -\stp_z &\mid \clF_{\stp_z }] < \infty.
\end{aligned}
\end{equation}
Then ${\bf \phi}$ is positive Harris recurrent.
\qed
\end{cor}

\section{Application to Stochastic Stabilization over an Erasure Channel}\label{sectionApplication}

\begin{figure}[h]
\centering
\Ebox{1.00}{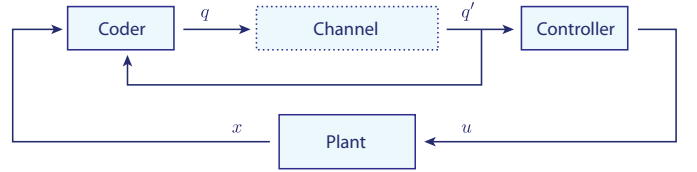}
\caption{Control over a discrete erasure channel with feedback. Coder represents the quantizer and the encoder. \label{LL}}
\end{figure}

The results of the previous section are now applied to a remote stabilization problem, in which the plant is open-loop unstable, and the controller has access to measurements from an erasure channel --- see Figure~\ref{LL}.

We begin with a scalar model;  extensions to the multivariate setting are contained in \Section{s:multi}.

\subsection{Scalar control/communication model}

Consider a scalar LTI discrete-time system described by
\begin{eqnarray}
\label{ProblemModel4}
x_{t+1}=ax_{t} + bu_{t} + d_t, \quad \quad t \geq 0
\end{eqnarray}
Here $x_t$ is the state at time $t$, $u_t$ is the control input, $x_0$ is a given initial condition, and $\{d_t \}$ is a sequence of zero-mean independent, identically distributed (i.i.d.) Gaussian random variables. 
It is assumed that     $|a| \geq 1$ and $b \neq 0$:  The system is open-loop unstable, but it is stabilizable.

This system is connected over an erasure channel with finite capacity to a controller, as shown in Figure~\ref{LL}. The controller has access to the information it has received through the channel. The controller estimates the state and then applies its control. We will establish bounds on data rates which lead to various versions of stochastic stability for the closed loop system.

The details of the communication channel are specified as follows:
The channel source consists of state values, taking values in $\Re$.  The source is \textit{quantized}:  The quantizer, at time $t \geq 0$, is represented by  a map $Q_t \colon\Re \to \Re$, characterized by a sequence of non-overlapping intervals ${\cal P}_t:= \{ {\cal B}_{i,t} \}$, with $|{\cal P}|=K+1$, such that $Q_t(x) = q^i_t$ if and only if $x \in \clB_{i,t}$;
that is, \[Q_t(x)=\sum_{i} q^i_t \times 1_{\{x \in \clB_{i,t}\}}.\]

The quantizer outputs are transmitted through a memoryless erasure channel, after being subjected to a bijective mapping, which is performed by the channel encoder: The channel encoder maps the quantizer output symbols to corresponding channel inputs $q \in \clM\eqdef\{1,2\dots,K+1\}$. An encoder at time $t$, denoted by ${\cal E}_t$, maps the quantizer outputs to $\clM$ such that ${\cal E}_t(Q_t(x_t))=q_t \in \clM$.

 The controller/decoder has access to noisy versions of the coder outputs for each time, which we denote by $\{q'\} \in \clM \cup \{e\}$, with $e$ denoting the erasure symbol, generated according to a probability distribution for every fixed $q \in \clM$.
 The channel transition probabilities are given by:
$$P(q'=i|q=i) = p, \quad \quad P(q'=e|q=i) = 1-p, \quad \quad i \in {\cal M}.$$

For each time $t \geq 0$, the controller/decoder applies a mapping ${\cal D}_t: \clM \cup \{e\} \to \mathbb{R}$, given by:
\[{\cal D}_t(q'_t) = {\cal E}_t^{-1}(q'_t) \times 1_{\{q'_t \neq e\}} + 0 \times 1_{\{q'_t=e\}}\]
We restrict the analysis to a class of uniform quantizers, defined by two parameters: bin size $\Delta>0$, and an even number $K\ge 2$.
The set $\clM$ consists of $K+1$ elements.
The uniform quantizer map  is defined as follows:  For $k=1,2\dots,K$,
\begin{eqnarray}
\!
Q_K^{\Delta}(x) = \begin{cases}  &  (k - \half (K + 1) ) \Delta,  \nonumber \\
&  \quad \quad \quad   \mbox{if} \ \ x \in [ (k-1-\half K   ) \Delta, (k-\half K  ) \Delta) \nonumber \\[.2cm]
&  (\half (K - 1) ) \Delta,   \quad \mbox{if} \ \ x = \half K \Delta \nonumber \\[.2cm]
&   0, \quad \quad \quad\quad  \quad\quad \mbox{if} \ \ x \not\in [- \half K  \Delta, \half K  \Delta].
\end{cases} \nonumber
\end{eqnarray}

We consider  quantizers that are \textit{adaptive},  so that the bin process can vary with time.
The   bin size  $\Delta_t$ at time $t$ is assumed to be a function of the previous value $\Delta_{t-1}$ and the past channel output $q'_{t-1}$.

\subsection{Stochastic stabilization over an erasure channel}

Consider the following time-invariant model. Let $\Upsilon_t$ denote an i.i.d. binary sequence of random variables, representing the erasure process in the channel, where the event $\Upsilon_t=1$ indicates that the signal is transmitted with no error through the erasure channel at time $t$. Let $p=\Expect[\Upsilon_t]$ denote the probability of success in transmission.

The following key assumptions are imposed throughout this section:  Given $K\ge 2$ introduced in the definition of the quantizer, define the \textit{rate variable}s
\begin{equation}\label{e:Rbdd}
R := \log_2(K+1) \quad \quad R'=\log_2(K)
\end{equation}
We fix positive scalars $\delta, \alpha $ satisfying $\alpha < 1$, and
\begin{equation}
\alpha > |a| 2^{-R'},
\label{e:Rbdd2}
\end{equation}
and
\begin{equation}
\alpha (|a|+\delta)^{p^{-1}-1} < 1.
\label{e:Rbdd3}
\end{equation}

We note that the (Shannon) capacity of such an erasure channel is given by $\log_2(K+1) p$ \cite{Cover}. From (\ref{e:Rbdd})-(\ref{e:Rbdd3}) it follows that if $\log_2(K) p > \log_2(|a|)$, then $\alpha, \delta$ exist such that the above are satisfied.


To define the bin-size update rule we require another constant $L>0$, chosen so that $L' =: \alpha L\ge 1$, where we take the lower bound as $1$ for convenience; any positive number would suffice.

Define the mapping $H: \mathbb{R} \times \mathbb{R} \times \{0,1\} \to \mathbb{R}$,
\begin{eqnarray}
 H(\Delta,h,p) = |a| + \delta \quad && \mbox{if } \quad |h| > 1, \quad \mbox{or} \quad \quad p = 0  \nonumber \\
 H(\Delta,h,p) = \alpha \quad && \mbox{if } \quad 0 \leq |h| \leq 1, p=1, \Delta \geq L   \nonumber \\
 H(\Delta,h,p) = 1 \quad \quad  && \mbox{if } \quad 0 \leq |h| \leq 1, p=1, \Delta < L \nonumber
\end{eqnarray}
Then with $\Delta_0 \geq L$ selected otherwise arbitrarily, define
\begin{eqnarray}
\label{QuantizerUpdate2}
u_t &=& - {a \over b} \hat{x}_{t}, \nonumber
\\
\hat{x}_t &=& {\cal D}_t(q'_t) = \Upsilon_t Q_K^{\Delta_t}(x_t), \nonumber
\\
\Delta_{t+1} &=& \Delta_t H(\Delta_t,|h_t|,\Upsilon_t), \text{\it where \ } h_t ={x_{t} \over \Delta_t 2^{R'-1}}
\end{eqnarray}
The update equations above imply that
\begin{eqnarray}
\label{definitionLprime}
\Delta_t \geq L \alpha =:L' \ge 1\, .
\end{eqnarray}

Given the channel output $q'_t \neq e$, the controller can deduce the realization of $\Upsilon_t$ and the event $\{|h_t| > 1\}$
simultaneously. This is due to the observation that if the channel output is not the erasure symbol, the controller knows that the signal is received with no error.  If $q'_t=e$, then the controller applies $0$ as its control input and enlarges the bin size of the quantizer.

\begin{lem} \label{MarkovChain}
Under (\ref{QuantizerUpdate2}), the process $(x_t,\Delta_t)$ is a Markov chain.
\end{lem}
\textbf{Proof:}
The system state dynamics can be expressed $ x_{t+1}= ax_t - a \hat{x}_t + d_t$,  where
$\hat{x}_{t}= \Upsilon_t Q_K^{\Delta_t} (x_{t}) $.  It follows that the pair process
$(x_t,\Delta_t)$ evolves as a nonlinear state space model,
\begin{equation}
\begin{aligned}
x_{t+1}  &= a (x_{t} - \Upsilon_t Q_K^{\Delta_t} (x_{t})) + d_t
\\[.2cm]
\Delta_{t+1}  &= \Delta_t H(\Delta_t, |{x_t \over 2^{R'-1}\Delta_t}|,\Upsilon_t)\, ,
\end{aligned}
\label{e:LemmaMCpf}
\end{equation}
in which $(d_t,\Upsilon_t)$ is i.i.d..
Thus, $(x_t, \Delta_t)$ form a Markov chain (see \cite[Ch.~2]{MCSS}).
\qed


Our result on the existence and uniqueness of an invariant probability measure is the following. 
\begin{thm}
\label{Inv}
For an adaptive quantizer satisfying \eqref{e:Rbdd}-\eqref{e:Rbdd3}, suppose that
the quantizer bin sizes are such that their base-2 logarithms are integer multiples of some scalar $s$,  and $\log_2(H( \varble ))$ takes values in integer multiples of $s$. Then the process $(x_t,\Delta_t)$ forms a positive Harris recurrent Markov chain, with a unique invariant probability measure $\pi$. If the integers taken are relatively prime (that is they share no common divisors except for $1$), then the invariant probability measure is independent of the value of the integer multiplying $s$.
\qed
\end{thm}

Under slightly stronger conditions we obtain a finite second moment:
\begin{thm}
\label{FiniteMoment}
Suppose that the assumptions of Theorem~\ref{Inv} hold, and  in addition we have the bound
\begin{eqnarray}\label{secondMomentBound0}
	a^2 \Bigl(1-p + {p \over (2^{R}-1)^2} \Bigr) < 1.
\end{eqnarray}
It then follows that for each initial condition $(x_0,\Delta_0)$,
$$
\lim_{t \to \infty} \Expect[x_t^2]  = \Expect_\pi[x_0^2] < \infty \, .
$$
\qed
\end{thm}

\begin{remark}
We note that Minero et al \cite{Minero}, in Theorem 4.1, observed that a necessary condition for mean square stability is that the following holds:
$$
	|a|^2 \Bigl(1-p + {p \over (2^{R})^2} \Bigr) < 1.
$$
Thus, our sufficiency proof almost meets this bound except for an additional transmitted symbol.
\end{remark}

We now consider the $m$th moment case. This moment can become useful for studying multi-dimensional systems for a sequential analysis of the modes, as we briefly discuss in \Section{s:multi}.

\begin{thm}
\label{FinitemMoment}
Consider the scalar system in (\ref{ProblemModel4}). Let $m \in \mathbb{N}$, suppose that the assumptions of Theorem~\ref{Inv} hold, and  in addition we have the inequality, \[|a|^m \Bigl(1-p + {p \over (2^{R}-1)^m} \Bigr) < 1.\]
It then follows that with the adaptive quantization policy considered, $\lim_{t \to \infty} \Expect[|x_t|^m]  = \Expect_\pi[|x_0|^m] < \infty \, .$
\qed
\end{thm}

\subsection{Connections with the drift criteria and the proof program}

Stability of the control/communication model is established using  the random-time stochastic drift criteria presented in the previous section, applied to   the Markov chain ${\bf \phi} = (x_t,\Delta_t)$ (see Lemma~\ref{MarkovChain}).   We provide an overview here, and the details can be found in the appendix.


\begin{figure}[h]
\centering
\Ebox{1.00}{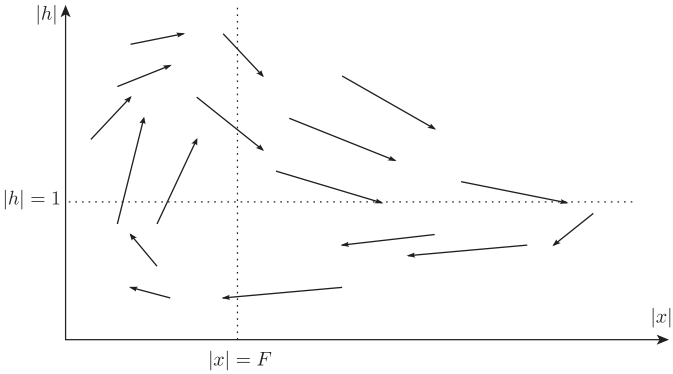}
\caption{Drift in the Markov Process. When under-zoomed, the error increases on average and the quantizer {\em zooms out}; when perfectly-zoomed, the error decreases and the quantizer {\em zooms in}.  \label{LLLLLLL}}
\end{figure}

Figure~\ref{LLLLLLL} provides some intuition on the construction of stopping times and the Lyapunov functions.
Recall that    $h_t= x_t / (2^{R'-1} \Delta_t)$
was introduced in  \eqref{QuantizerUpdate2}.
The arrows shown in the figure denote the mean one-step increments of $(x_t, h_t)$:  That is, the arrow $\nu$ with base at $(x,h)$ is defined by,
\[
\nu = \Expect[(x_{t+1}, h_{t+1})-(x_t, h_t)\mid (x_t, h_t) = (x,h)]
\]

With $F>0$ fixed,
and with $F'= F  2^{-(R'-1)}$,
two sets are used to define the small set in the drift criteria,
$C_x= \{x: |x| \leq F\}$ for $F>0$,  and  $C_{\Delta} = \{\Delta: \Delta \leq F'\}$.
Denote  $C_h=\{h: |h| \leq 1\}$, and assume that $F>0$ is chosen sufficiently large so that
$(x_t,\Delta_t) \in C_x \times C_{\Delta}$ whenever $(x_t,h_t) \in C_{x} \times C_h$.
When $x_t$ is outside $C_x$ and $h_t$ outside $C_h$ (the {\em under-zoomed phase} of the quantizer), there is a drift for $h_t$ towards $C_h$.   When the process $x_t$ reaches $C_h$ (the {\em perfectly-zoomed phase} of the quantizer), then the process
drifts towards $C_x$.

%
%
%
%

We next construct the sequence of stopping times required in the drift criteria of Section~\ref{driftCriteria}.  The controller can receive meaningful information regarding the state of the system when two events occur concurrently: the channel carries information with no error, and the source lies in the granular region of the quantizer: That is,   $x_t \in [- \half K  \Delta_t, \half K  \Delta_t]$ (or $|h_t| \leq 1$) and $\Upsilon_t=1$. The stopping times are taken to be the times at which both of these events occur. With,
$|h_0| \leq 1, \Upsilon_0=1$, we define $\stp_0 = 0$ and
\begin{eqnarray}
  \stp_{z+1} &=& \inf \{k > \stp_z : |h_{k}| \leq  1, \Upsilon_k=1 \}, \quad z \in \mathbb{N}. \nonumber 
\end{eqnarray}

These stopping times are nearly geometric when the bin size is large.  The proof of Proposition \ref{stoppingTimeDistribution} is presented in \Section{ProofStocStabSystemTimeDistribution}.

\begin{prop}
\label{stoppingTimeDistribution}
The discrete probability measure $\Prob(\stp_{z+1}-\stp_z=k\mid x_{\stp_z },\Delta_{\stp_z })$ satisfies,
\[
(1-p)^{k-1} \leq \Prob(\stp_{z+1}-\stp_z \geq k | x_{\stp_z },\Delta_{\stp_z }) \leq (1-p)^{k-1} +  o(1),
\]
where $o(1) \to 0$ as $\Delta_{\stp_z} \to \infty$ uniformly in $x_{\stp_z}$. \qed
\end{prop}

The next step is to establish irreducibility structure.
The proof of the following is contained in \Section{ProofStocStabSystemSmallSet}.

\begin{prop} \label{compactSetSmall}
Under the assumptions of \Theorem{Inv},  the chain $(x_t,\Delta_t)$ is
$\varphi$-irreducible for some $\varphi$,  it is  aperiodic, and all compact sets are small.
\qed
\end{prop}

We now provide sketches of the proofs of the main results.  The details are collected together in the appendix.

\textit{\textbf{Sketch of Proof of \Theorem{Inv}:}}
The logarithmic function $V_0(x_t,\Delta_t)=  \log(\Delta^2) + B_0$ for some $B_0 > 0$ serves as the Lyapunov function in  \eqref{e:thm5delta}, with $f(x,\Delta)$ set as a constant. Note that by (\ref{definitionLprime}), $V_0(x_t,\Delta_t) > 0$.

Together with Propositions~\ref{stoppingTimeDistribution} and \ref{compactSetSmall}, we apply Theorem~\ref{thm5} in the special case of Corollary~\ref{corol}. Proposition~\ref{compactSetSmall} implies the existence of a unique invariant measure. Details of the proof are presented in section \ref{ProofStocStabSystemProofInv}. \qed

\textit{\textbf{Sketch of Proof of \Theorem{FiniteMoment}:}} A quadratic Lyapunov function $V_2(x_t,\Delta_t)= \Delta^2_t$ is used, along with
 $\delta(x_t,\Delta_t)=\epsilon \Delta_t^2$ for some $\epsilon>0$, and  $f(x_t,\Delta_t) = \xi x_t^2$ for some $\xi > 0$.
The bound \eqref{e:thm5delta} is established in \Section{ProofStocStabSystemFiniteMoment},
so that the limit $\lim_{t \to \infty} \Expect[x_t^2] $ exists and is finite by Theorem~\ref{thm5}. Details of the proof are presented in section \ref{ProofStocStabSystemFiniteMoment}.
 \qed

\textit{\textbf{Sketch of Proof of \Theorem{FinitemMoment}:}} Theorem~\ref{thm5} is applied with the Lyapunov function $V_m(x_t,\Delta_t)= \Delta^m_t$, $\delta(x_t,\Delta_t)=\epsilon \Delta_t^m$ for some $\epsilon>0$, and  $f(x_t,\Delta_t) = \xi |x_t|^m$ for some $\xi > 0$. Details of the proof are presented in section \ref{ProofStocStabSystemFinitemMoment}.
 \qed

\subsection{Simulation}

\begin{figure}[h]
\Ebox{1.00}{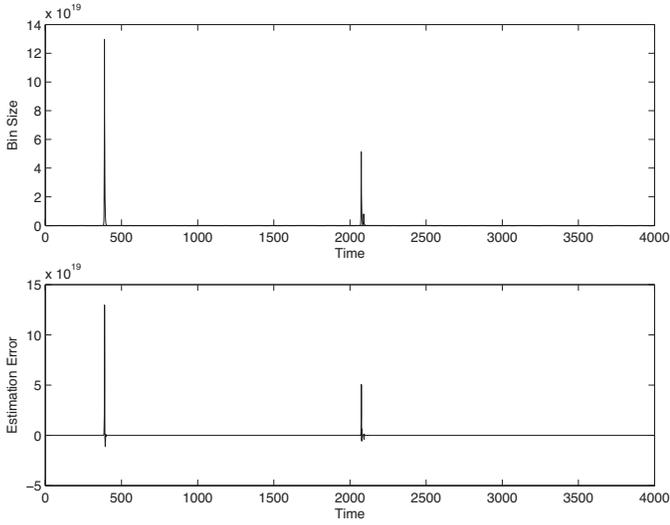}
\caption{Sample path for a stochastically stable system with a 5-bin quantizer. \label{LLL}}
\end{figure}

\begin{figure}[h]
\Ebox{1.00}{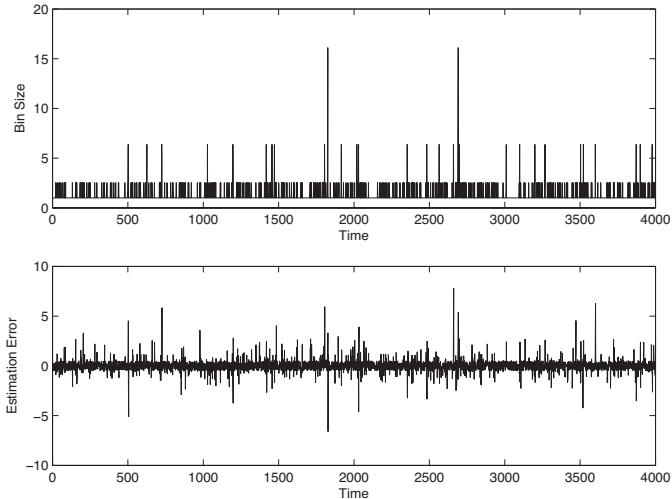}
\caption{Sample path with a 17-bin quantizer; a more desirable path. \label{LLLL}}
\end{figure}

Consider a linear system
$$
x_{t+1}= a x_t + u_t + d_t,
$$
with $a=2.5$, $\{d_t\}$ is an i.i.d.\ $N(0,1)$  Gaussian sequence. The erasure channel has erasure probability $1-p=0.1$.
For stability with a finite second moment, we employ a quantizer with rate
\[
\log_2( \lceil \sqrt{p \over {1 \over a^2} - (1-p)} \rceil +1)=\log_2(5)
\]
bits. That is, a uniform quantizer with $5$ bins. We have taken $L'=1$.
  Figures~\ref{LLL} and \ref{LLLL} illustrate the conclusions of the stochastic stability results
presented in Theorems~\ref{Inv} and \ref{FiniteMoment}.  The plots show  the under-zoomed and perfectly-zoomed phases, with the peaks in the plots showing the under-zoom phases. For the plot with $5$ levels, the system is positive Harris recurrent, since the update equations are such that $\alpha=0.629$, $\delta=0.025$ and $\log_2(H(\cdot)) \in \{-0.6744, 0, 1.363\}$. These values satisfy the irreducibility condition since $-0.6744 = - 1.363/2$,  and hence the communication conditions in Theorem~\ref{Inv} are satisfied. Furthermore, \eqref{e:Rbdd3} is satisfied since
\[\alpha ({|a|+\delta})^{p^{-1}-1} <  0.698 < 1.
\]
Increasing the bit rate by only two bits in Figure~\ref{LLLL} leads to a much more desirable sample path, for which the severity of rare events is reduced.



%

\subsection{Extension to multi-dimensional systems}
\label{s:multi}

The control laws and analysis can be extended to  multi-dimensional models. Consider the multi-dimensional linear system
\begin{eqnarray}\label{vectorEqn1}
x_{t+1}=Ax_t + Bu_t + d_t,
\end{eqnarray}
where $x_t \in \mathbb{R}^n$ is the state at time $t$, $u_t$ is the control input, and $\{d_t \}$ is a sequence of zero-mean independent, identically distributed (i.i.d.) $\mathbb{R}^n$-valued zero-mean Gaussian random variables. Here $A$ is the system matrix with at least one eigenvalue greater than $1$ in magnitude, so that the system is open-loop unstable.  Without any loss of generality, we assume $A$ to be in Jordan form.  Suppose that $B$ is invertible for ease in presentation.

The approach for the scalar systems is applicable, however some extension is needed. Toward this goal, one can adopt two approaches.

In one scheme, one could consider a sequential stabilization of the scalar components. In particular, we can perform an analysis by considering a lower to upper sequential approach, considering stabilized modes at particular stopping times as noise with a finite moment. Using an inductive argument, one can first start with the lowest mode (in the matrix diagonal) of the system, and stabilize that mode so that there is a finite invariant $m$th moment of the state. We note that the random process for the upper mode might not have Markov dynamics for its marginal, but the joint system consisting of all modes and quantizer parameters will be Markov. Such a sequential scheme ensures a successful application of the scalar analysis presented in this paper to the vector case. One technicality that arises in this case is the fact that the effective disturbance affecting the stochastic evolution of a repeated mode in a Jordan block is no longer Gaussian, but can be guaranteed to have a sufficiently light tail distribution by Theorem~\ref{FinitemMoment}.

Another approach is to adopt the discussions in \cite{YukTutorial2011} and Section IV of \cite{YukITA2011} (see also 
\cite{Matveev} for related constructions) and apply a vector quantizer by transmitting the quantizer bits for the entire $\mathbb{R}^n$-valued state. In particular, by defining a vector quantizer as a product of scalar quantizers along each (possibly generalized) eigenvector with a common under-zoom bin, letting $h^i$ denote the ratio of the state and the bin range of the corresponding $i$th scalar quantizer and defining a sequence of stopping times as follows with $\stp_0 = 0$ and for $z \in \mathbb{N}$
\begin{eqnarray}
  \stp_{z+1} = \inf \{k > \stp_z : |h^i_{k}| \leq  1, i=1,2\dots,n, \quad \Upsilon_k=1 \}, \nonumber 
\end{eqnarray}
the analysis can be carried over through a geometric bound on the distribution of subsequent stopping times, obtained by an application of the union bound. See \cite{YukITA2011} for details.

\section{Concluding Remarks}
\label{s:conc}

This paper contains two main contributions. One is on a general drift approach for verifying stochastic stability of Markov Chains. The other is on stabilization over erasure channels. We believe that the results presented in this paper will have many applications within the context of network stability, and networked control systems as well as information theoretic applications. Important previous research on performance bounds for variable-length decoding schemes use stopping time arguments \cite{Burnashev}, \cite{PolyanskiyISIT10}, and this could form the starting point of a Lyapunov analysis. 

The methods of this paper can be extended to a large class of networked control systems and delay-sensitive information transmission. For networked control systems, the effects of randomness in delay for transmission of sensor or controller signals (see for example \cite{Montestruque}, \cite{Cloosterman}) is an application area where the research reported in this paper can be applied.


Rates of convergence under random-time drift is one direction of future research. It is apparent that the nature of the drift as well as the distribution of stopping times used for drift will play a role in the rate of convergence. We refer the reader to \cite{Douc} and \cite{ConnorFort}, for results when the drift times are deterministic.

Positive Harris recurrence can be a crude measure of stability, as seen in the numerical results in this paper. Computing the sensitivity of performance to the bit rate is an important future research problem for practical applications. For example, it was observed in Figure~\ref{LLLL}  that increasing the bit rate by only two bits leads to   much more desirable sample paths, and the magnitudes of rare events are significantly reduced.   A Markovian framework is  valuable for sensitivity analysis, as applied in the reinforcement learning   literature (see commentary in Section~17.7 of \cite{MCSS}).

%


\section{Appendix}

\subsection{Proofs of the Stochastic Stability Theorems}
\label{ProofStocStab}

\subsubsection{Proof of Theorem~\ref{thm5} (i)}

The proof is similar to the proof of the Comparison Theorem of \cite{MCSS}:
Define the sequence $\{M_z : z\ge 0\}$ by $M_{0} = V(\phi_{0})$, and for $z\ge0$,
\[
M_{z+1} = V(\phi_{\stp_{z+1}}) + \sum_{k=0}^{z} (\delta(\phi_{\stp_k}) - b1_{\{\phi_{\stp_k} \in C\}} )
\]
Under the assumed drift condition we have,
\[
\Expect[ M_{z+1} \mid \clF_{\stp_z}]
	\leq  V(\phi_{\stp_z }) + \sum_{k=0}^{z-1} (\delta(\phi_{\stp_k}) - b1_{\{\phi_{\stp_k} \in C\}}),
\]
which implies the super-martingale bound,
\[
\Expect[ M_{z+1} \mid \clF_{\stp_z}] \leq M_z
\]


For a measurable subset $C\subset\state$ we denote the first hitting time for the sampled chain,
\begin{equation}
\tn_C= \min\{ z \ge 1 : \phi_{\stp_z} \in C \}
\label{e:hatauC}
\end{equation}
Define   $\tn_C^n  = \min (n, \tn_C)$ for any $n\ge 1$.
Then  $\Expect[M_{\tn_C^n}] \leq M_0$ for any $n \in \mathbb{Z}$,
 and
\[
\Expect[\sum_{k=0}^{\tn_C^n-1} \delta(\phi_{\stp_k}) |\clF_0 ] \leq M_0 +b.
\]
Applying the bound
$\Expect \Bigl[\sum_{k=\stp_z}^{\stp_{z+1}-1} f(\phi_k)  \mid \clF_{\stp_z }\Bigr]  \le \delta(\phi_{\stp_z})$ and that $f(\phi) \geq 1$,
the following bound is obtained from the smoothing property of the conditional expectation:
\[
\begin{aligned}
 \Expect[\stp_{\tn_C^n} \mid \clF_0 ]
	&=  \Expect\Bigl[\sum_{i=0}^{\tn_C^n-1} \Expect[\stp_{i+1} - \stp_i] |\clF_0 ] \Bigr]
	\\
	& \leq \Expect\Bigl[\sum_{i=0}^{\tn_C^n-1} \delta(\phi_{\stp_i}) |\clF_0 \Bigr] \leq M_0 +b
\end{aligned}
\]
Hence by the monotone convergence theorem,
\[
\Expect[\tau_C]
\le
\Expect[\stp_{\tn_C}] =\lim_{n \to \infty} \Expect[\stp_{\tn_C^n} \mid \clF_0 ] \leq M_0 +b.
\]
Consequently we obtain that
\[\sup_{\phi \in C} \Expect[\tau_C] < \infty, \]
as well as recurrence of the chain, $\Prob_{\phi}(\tau_C < \infty ) =1 $ for any $\phi\in\state$.
Positive Harris recurrence  now follows from \cite{Meyn} Thm. 4.1.
\qed

The following result is key to obtaining moment bounds. The inequality \eqref{e:V3} is known as \textit{drift condition (V3)} \cite{MCSS}. Define,
\begin{equation}
V^*_f(\phi) \eqdef\Expect_{\phi}\Bigl[ \sum_{t=0}^{\tau_C-1} f(\phi_t)\Bigr]\, \quad \phi\in\state.
\label{e:Vf}
\end{equation}
\begin{lem}
\label{l:thm5}
Suppose that ${\bf \phi}$ satisfies all of the assumptions of Theorem~\ref{thm5}, except that the $\psi$-irreducibility assumption is relaxed.  Then, there is a constant $b_f$ such that the following bounds hold
\begin{eqnarray}
PV^*_f &\le& V^*_f -f + b_f\ind_C
\label{e:V3}
\\
V^*_f(\phi) &  \le & V(\phi) + b_f,\qquad \phi\in\state.
\label{e:VfBdd}
\end{eqnarray}
\end{lem}

\subsubsection*{Proof}
The drift condition \eqref{e:V3} is given in Theorem~14.0.1 of \cite{MCSS}.

The proof of \eqref{e:VfBdd} is based on familiar super-martingale arguments:   Denote $M_{0} = V(\phi_{0})$, and
and for $z\ge 0$,
\begin{eqnarray}
\label{boundMomentProof1}
M_{z+1}
= V(\phi_{\stp_{z+1}})
	- \sum_{k=0}^{\stp_{z+1}-1} \bigl( - f(\phi_k)  + b1_{\{\phi_{\stp_k} \in C\}} \bigr)
\end{eqnarray}
The super-martingale property for $\{ M_z\}$   follows from the assumed drift condition:
\begin{eqnarray}\label{boundMomentProof2}
&& \Expect[ M_{z+1} \mid \clF_{\stp_z}] = M_z + \Expect\Bigl[   V(\phi_{\stp_{z+1} }) -  V(\phi_{\stp_z }) \nonumber \\
&& \quad \quad \quad + \sum_{k=\stp_z}^{\stp_{z+1}-1}  (f(\phi_k) - b1_{\{\phi_{\stp_k} \in C\}}) \mid \clF_{\stp_z} \Bigr]
\leq M_z
\end{eqnarray}
As in the previous proof we bound expectations involving the stopping time $\tn_C$ beginning with its truncation $\tn_C^n  = \min (n, \tn_C)$.

The super-martingale property gives $\Expect[M_{\tn_C^n}] \leq M_0$,  and once again it
 follows again by the monotone convergence theorem that $V_f^*$ satisfies the bound \eqref{e:VfBdd} as claimed.  
\qed

\subsubsection{Proof of Theorem~\ref{thm5} (ii) and (iii)}   The existence of a finite moment follows from Lemma~\ref{l:thm5} and the following generalization of Kac's Theorem (see \cite[Theorem 10.4.9]{MCSS}):
\begin{equation}
\pi(f):= \int \pi(d\phi) f(\phi) = \int_{A} \pi(d\phi) \Expect_{\phi}\bigg( \sum_{t=0}^{\tau_A-1} f(\phi_t)\bigg),
\label{e:KAC}
\end{equation}
where $A$ is any set satisfying $\pi(A)>0$,  and $\tau_A = \inf(t \ge 1: \phi_t \in A)$. The super-martingale argument above ensures that the expectation under the invariant probability measure is bounded by recognizing $C$ as a recurrent set.

(iii) now follows from the ergodic theorem for Markov chains, see \cite[Theorem 17.0.1]{MCSS}.
\qed

\subsubsection{Proof of Theorem~\ref{thm25Feller}}

The existence of an invariant probability measure in (i) follows from Theorem 12.3.4 of \cite{MCSS} (the solution to the drift condition (V2) can be taken to be the mean hitting time,  $V(\phi) = \Expect_{\phi}[\tau_C]$).  See also \cite[Theorem 3.1]{Lasserre}.

The proof of (ii) is similar.  Rather than work with the mean return time to $C$, we consider the function $V_f^*$ defined in Lemma~\ref{l:thm5}.  We have by the Comparison Theorem of \cite{MCSS},
\[
0\le P^n V^*_f \le V^*_f + nb_f - \sum_{t=0}^{n-1} P^t f
\]
Hence for any $\phi\in\state$,
\begin{equation}
\limsup_{n\to\infty} \frac{1}{n}  \sum_{t=0}^{n-1} P^t f \, (\phi) \le b_f.
\label{e:V3f_bdd}
\end{equation}
Suppose that $\pi$ is any invariant probability measure.  Fix $N<\infty$, let $f_N=\min(N,f)$,  and apply Fatou's Lemma as follows,
\[
\begin{aligned}
\pi(f_N) & = \limsup_{n\to\infty} \pi\Bigl( \frac{1}{n}  \sum_{t=0}^{n-1} P^t f_N \Bigr)
\\
& \le\pi\Bigl(  \limsup_{n\to\infty} \frac{1}{n}  \sum_{t=0}^{n-1} P^t f_N \Bigr) \le b_f\, .
\end{aligned}
\]
Fatou's Lemma is justified to obtain the  first inequality, because $f_N$ is bounded. The second inequality holds by \eqref{e:V3f_bdd} and since $f_N\le f$.   The monotone convergence theorem then gives $\pi(f)\le b_f$.
 \qed

\subsection{Proofs of  Stability: Stochastic Stabilization over an Erasure Channel }
\label{ProofStocStabSystem}

\subsubsection{Proof of Proposition \ref{stoppingTimeDistribution}} \label{ProofStocStabSystemTimeDistribution}
We obtain upper and lower bounds below.

\begin{lem}
\label{KeyBound}
The discrete probability measure $\Prob(\stp_{z+1}-\stp_z=k\mid x_{\stp_z },\Delta_{\stp_z })$ satisfies
\[
\Prob(\stp_{z+1}-\stp_z \geq k | x_{\stp_z },\Delta_{\stp_z }) \leq (1-p)^{k-1} +  G_k(\Delta_{\stp_z}),
\]
where $G_k(\Delta_{\stp_z}) \to 0$ as $\Delta_{\stp_z} \to \infty$ uniformly in $x_{\stp_z}$.
\qed
\end{lem}

\textbf{Proof:}
Denote for $k \in \mathbb{N}$,
\begin{eqnarray}
&&\Theta_k := \Prob(\stp_{z+1} - \stp_z \geq k\mid x_{\stp_z },\Delta_{\stp_z }) \nonumber \\
&& \quad \quad = \Prob_{x_{\stp_z },\Delta_{\stp_z }}(\stp_{z+1} - \stp_z \geq k).
\label{e:stpP}
\end{eqnarray}
Without any loss, let $z=0, \stp_0=0$, so that $\Theta_k =\Prob_{x_0,\Delta_{0}}(\stp_{1} \geq k) $.

Now, at time $0$, upon receiving a message successfully, the estimation error satisfies $|x_0- \hat{x}_0| \leq \Delta_0/2$, as such we have that $|a||x_0 + (b/a)u_0| \leq |a|\Delta_0/2$.

The probability $\Theta_k$ for $k \geq 2$ is bounded as follows:
\begin{eqnarray}
&& \Theta_k = \Prob_{x_0,\Delta_{0}}\bigg( \bigcap_{s=1}^{k-1} (\Upsilon_s = 0) \cup (|h_s| > 1) \bigg) \nonumber \\
&\leq&\Prob_{x_0,\Delta_{0}}\bigg( \bigcap_{s=1}^{k-1} (\Upsilon_s = 0) \cup (|x_s| \geq 2^{R'-1} (|a|+\delta)^{s-1} \alpha \Delta_0) \bigg) \nonumber \\
&=& \Prob_{x_0,\Delta_{0}}\bigg( \bigcap_{s=1}^{k-1} (\Upsilon_s = 0) \nonumber \\
&&  \quad \quad \cup (|a^{s}(x_0 + \sum_{i=0}^{s-1}a^{-i-1}d_i)| \geq (|a|+\delta)^{s-1}2^{R'-1} \alpha\Delta_0) \bigg)  \nonumber \\
& \leq& \Prob_{x_0,\Delta_{0}}\bigg( \bigcap_{s=1}^{k-2} (\Upsilon_s = 0) \cup (|h_s| > 1) \,\Big|\,   \Upsilon_{k-1}=0 \bigg) (1- p)
\nonumber \\
&& + \Prob_{x_0,\Delta_{0}}\bigg( \bigcap_{s=1}^{k-2} (\Upsilon_s = 0) \cup (|h_s| > 1) \, \nonumber \\
&& \quad \quad \quad \quad \quad \Big|\,   |a^{k-1}(x_0 + (b/a)u_0 + \sum_{i=0}^{k-2}a^{-i-1}d_i)| \nonumber \\
&& \quad \quad \quad \quad \quad \quad \quad \quad \geq (|a|+\delta)^{k-2}2^{R'-1} \alpha \Delta_0  \bigg)
\nonumber \\
&& \quad {} \times \Prob_{x_0,\Delta_{0}}\bigg( |a^{k-1}(x_0 + (b/a)u_0 + \sum_{i=0}^{k-2}a^{-i-1}d_i)| \nonumber \\
&& \quad \quad \quad \quad \quad \quad \quad \quad \geq (|a|+\delta)^{k-2}2^{R'-1} \alpha \Delta_0 \bigg) \label{intUnionBound} \\
&\le&  \Prob_{x_0,\Delta_{0}} \bigg( \bigcap_{s=1}^{k-2} (\Upsilon_s = 0) \cup (|h_s| > 1) \mid \Upsilon_{k-1}=0 \bigg) (1- p) \nonumber \\
&& {}  + \Prob_{x_0,\Delta_{0}} \bigg( |a^{k-1}(x_0 + (b/a)u_0 + \sum_{i=0}^{k-2}a^{-i-1}d_i)| \nonumber \\
&& \quad \quad\quad \quad \quad \quad \quad \quad \geq (|a|+\delta)^{k-2}2^{R'-1} \alpha \Delta_0 \bigg) \nonumber \\
&=&  \Prob_{x_0,\Delta_{0}}(\stp_{1}  \geq k-1) (1- p)  \nonumber \\
&&{} + \Prob_{x_0,\Delta_{0}} \bigg( |a^{k-1}(x_0 + (b/a)u_0 + \sum_{i=0}^{k-2}a^{-i-1}d_i)| \nonumber \\
&& \quad \quad\quad \quad \quad \quad \quad \quad \geq (|a|+\delta)^{k-2}2^{R'-1} \alpha \Delta_0 \bigg).
\end{eqnarray}
In the above derivation, (\ref{intUnionBound}) follows from the following: For any three events $M,C,D$ in a common probability space
\[
\Prob\big( M \cap (C \cup D) \big) = \Prob\big( (M \cap C) \cup  (M \cap D) \big)
\le
 \Prob\big( M \cap C  \big)  +  \Prob\big(  M \cap D \big)
 \]
Now, observe that for $k\geq 2$,
\begin{eqnarray}
&& \Prob_{x_0,\Delta_{0}}\bigg(|(x_0 + (b/a)u_0 + \sum_{i=0}^{k-2}a^{-i-1}d_i)| \nonumber \\
&& \quad \quad \quad \quad \quad \quad \geq ({|a|+\delta \over |a|})^{k-2}2^{R'-1} {\alpha \over |a|} \Delta_0 \bigg) \nonumber \\
&& \le 2\Prob_{x_0,\Delta_{0}}\bigg(\sum_{i=0}^{k-2}a^{-i-1}d_i \nonumber \\
&& \quad \quad \quad \quad  \quad \quad  \quad \quad  \geq  (2^{R'-1}({|a|+\delta \over |a|})^{k-2} {\alpha \over |a|} - {1 \over 2} )\Delta_0\bigg) \nonumber \\
&& \le C  \exp\Bigl( -{((\xi^{k-2}N - 1/2)\Delta_0)^2 \over 2\sigma'^2} \Bigr)
\label{comp},
\end{eqnarray}
where (\ref{comp}) follows from (\ref{e:Rbdd2}), for this condition ensures that the term \[(2^{R'-1}({|a|+\delta \over |a|})^{k-2} {\alpha \over |a|} - {1 \over 2} )\]
 is positive for $k \geq 2$, and
 bounding the complementary error function by the following:
$\int_q^{\infty} \mu(dx) \leq  q^{-1} \int_q^{\infty}  x  \mu(dx)$,
for $q>0$.
In the above derivation, the constants are:
\[ \sigma'^2 = {E[d_1^2]\over 1 - |a|^{-2}}, \quad \xi = {|a|+\delta \over |a|}, \quad N = {2^{R'-1} \over (|a|)/\alpha },\] and
\[ C = 2\sigma' {1 \over \sqrt{2 \pi} (2N-1)\Delta_0/2 }.\]
Let us define:
\[\Xi_k :=  {((\xi^{k-2}N - 1/2)\Delta_0)^2 \over 2\sigma'^2}\]
and
\[\widetilde{\Xi}_k :=  {((\xi^{k}N - 1/2)\Delta_0)^2 \over 2\sigma'^2}\]

We can bound the probability $\Theta_k$ defined in \eqref{e:stpP}.  Since a decaying exponential decays faster than any decaying polynomial, for each $m \in \mathbb{N}_+$, there exists an $M< \infty$ such that for all $k \in \mathbb{N}$,
\begin{eqnarray}
C  e^{-\Xi_k} \leq M \widetilde{\Xi}_k^{-m}.
\end{eqnarray}
Thus, we have that
\begin{eqnarray}\label{expPoly}
&& \Prob_{x_0,\Delta_{0}}\bigg(x_0 + \sum_{i=0}^{k-2}a^{-i-1}d_i \geq ({|a|+\delta \over |a|})^{k-2}2^{R'-1} {\alpha \over |a|} \Delta_0 \bigg) \nonumber \\
&& \quad \quad \leq M \widetilde{\Xi}_k^{-m}.
\end{eqnarray}
Now $\Theta_1=1$
by definition, and
for $k >1 $,
\begin{eqnarray}\label{MultiBound1}
\Theta_k \leq \Theta_{k-1}(1- p) + C  e^{-\Xi_k}.
\end{eqnarray}
We obtain,
\begin{eqnarray}\label{MultiBound2}
\Theta_{k} &\le& \Theta_1 (1-p)^{k-1} + \sum_{s=1}^{k-1} (1-p)^{k-s-1} C  e^{-\Xi_s} \nonumber \\
&\le& \Theta_1 (1-p)^{k-1} +  \sum_{s=1}^{k-1} M (1-p)^{k-s-1} \widetilde{\Xi}_s^{-m} \nonumber \\
&=&  (1-p)^{k-1} +  G_k(\Delta_{\stp_0}),
\end{eqnarray}
where
\begin{eqnarray}\label{definitionofG}
G_k(\Delta_{\stp_0}) := \sum_{s=1}^{k-1} M (1-p)^{k-s-1} \widetilde{\Xi}_s^{-m}
\end{eqnarray}
It now follows that,
\begin{eqnarray}
&& G_k(\Delta_{\stp_0}) = \sum_{s=1}^{k-1} M (1-p)^{k-s-1} \widetilde{\Xi}_s^{-m} \nonumber \\
&& = \Delta_{0}^{-2m} \sum_{s=1}^{k-1} M (1-p)^{k-s-1} \bigg({(\xi^{s}N - 1/2)^2 \over (2\sigma'^2)}\bigg)^{-m} \nonumber \\
&& = \Delta_{0}^{-2m} (1-p)^{k-1} \sum_{s=1}^{k-1} M (1-p)^{-s} \bigg({(\xi^{s}N - 1/2)^2 \over (2\sigma'^2)} \bigg)^{-m} \nonumber \\
&& = \Delta_{0}^{-2m} (1-p)^{k-1} \nonumber \\
&& \quad \quad \quad \times \sum_{s=1}^{k-1} M (1-p)^{-s} \xi^{-2ms} (N - {1 \over 2\xi^{s}})^{-2m} (2\sigma'^2)^{m} \nonumber \\
&& \leq \Delta_{0}^{-2m} (1-p)^{k-1} \Gamma_m \sum_{s=1}^{k-1}  \bigg((1-p)\xi^{2m}\bigg)^{-s} \nonumber \\
&& \leq \Gamma_m  \Delta_{0}^{-2m} (1-p)^{k-1} {\bigg((1-p)\xi^{2m}\bigg)^{-k} - 1  \over \bigg((1-p) \xi^{2m} \bigg)^{-1} - 1},
\end{eqnarray}
with $\Gamma_m = M(N - {1 \over 2\xi})^{-2m}(2\sigma'^2)^{m} < \infty$. It follows that if $m$ is taken such that
\begin{eqnarray}\label{conditionOnMomentBoundedness2}
(1-p)\xi^{2m} > 1,
\end{eqnarray}
then $\lim_{\Delta_{0} \to \infty} G_k(\Delta_{0}) = 0$, and for all $k \in \mathbb{N}$
\begin{eqnarray}\label{conditionOnMomentBoundedness}
\Theta_k \leq (1-p)^{k-1} \bigg(1 +  \Gamma_m \Delta_{\stp_0}^{-2m} {1 \over 1 - ((1-p)\xi^{2m})^{-1} }\bigg)
\end{eqnarray}
\qed

\begin{lem}\label{KeyBound2}
The discrete probability measure $\Prob(\stp_{z+1}-\stp_z=k\mid x_{\stp_z },\Delta_{\stp_z })$ satisfies
\[\Prob(\stp_{z+1}-\stp_z \geq k | x_{\stp_z },\Delta_{\stp_z }) \geq (1-p)^{k-1},\]
for all realizations of $x_{\stp_z },\Delta_{\stp_z }$.
\qed
\end{lem}

\textbf{Proof:}
This follows since
$$\Prob_{x_0,\Delta_{0}}\bigg( \bigcap_{s=1}^{k-1} (\Upsilon_s = 0) \cup (|h_s| > 1) \bigg) \geq \Prob_{x_0,\Delta_{0}}\bigg( \bigcap_{s=1}^{k-1} (\Upsilon_s = 0)\bigg).$$
\qed
As a consequence of Lemma~\ref{KeyBound} and Lemma~\ref{KeyBound2}, the probability $\Prob( \stp_{z+1}-\stp_z = k \mid x_{\stp_z},\Delta_{\stp_z})$
tends to $ (1-p)^{k-1} p$ as  $\Delta_{\stp_z} \to \infty$. \qed

\subsubsection{Proof of Proposition \ref{compactSetSmall}}\label{ProofStocStabSystemSmallSet}

Let the values taken by $\log_2(H(\cdot))/s$ be $\{-\tilde{A}, 0, \tilde{B}\}$. Let
\begin{eqnarray}
&& \mathbb{L}_{z_0,\tilde{A},\tilde{B}} := \{n \in \mathbb{N}, n \geq \log_2({L'})/s \nonumber \\
&& \quad \quad : \exists N_{\tilde{A}} \in \mathbb{N}, N_{\tilde{B}} \in \mathbb{N}, n
	= -N_{\tilde{A}} \tilde{A} + N_{\tilde{B}} \tilde{B} + z_0 \}.\nonumber
\end{eqnarray}
Since we have by (\ref{QuantizerUpdate2}) \[\Delta_{t+1}=H(\Delta_t,|h_t|,\Upsilon_t) \Delta_t,\] it follows that,
$$
 	\log_2(\Delta_{t+1})/s = \log_2(H(\Delta_t,|h_t|,\Upsilon_t))/s + \log_2(\Delta_{t})/s
 $$
 is also an integer. We will establish that $\mathbb{L}_{z_0,\tilde{A},\tilde{B}}$ forms a communication class, where $z_0 = \log_2(\Delta_{0})/s$ is the initial condition of the parameter for the quantizer.
Furthermore, since the source process $x_t$ is ``Lebesgue-irreducible'' (for the system noise admits a probability density function that is positive everywhere) and there is a uniform lower bound $L'$ on bin-sizes, the error process takes values in any of the admissible quantizer bins with non-zero probability. In view of these, we now establish that the Markov chain is irreducible.

Given $l,k \in \mathbb{L}_{z_0,\tilde{A},\tilde{B}}$, there exist $N_{\tilde{A}}, N_{\tilde{B}} \in \mathbb{N}$ such that $l-k = - N_{\tilde{A}} \tilde{A} + N_{\tilde{B}} \tilde{B}$.
In particular, if at time $0$, the quantizer is perfectly zoomed and $\Delta_0 = 2^{sk}$, then there exists a sequence of events consisting of $N_{\tilde{B}}$ erasure events (simultaneously satisfying $|h_t| \leq 1$) and consequently $N_{\tilde{A}}$ zoom-in events taking place with probability at least:
\begin{eqnarray}\label{MoveBin2Bin}
&& \bigg( p \Prob(d_t \in [-(\alpha 2^{R'} - |a|)L'/2, (\alpha 2^{R'} - |a|)L'/2]) \bigg)^{N_{\tilde{A}}} \nonumber \\
&& \quad \quad \quad \quad \times \bigg( (1-p) \Prob(|d_t| \leq \delta 2^{R'-1} L') \bigg)^{N_{\tilde{B}}} \nonumber \\
&& > 0,
 \end{eqnarray}
so that $\Prob(\Delta_{N_{\tilde{A}} + N_{\tilde{B}}} = 2^{sl} | \Delta_0=2^{ks},x_0) > 0$, uniformly in $x_0$. In the following we will consider this sequence of events.

Now, for some distribution ${\cal K} $ on positive integers, $E \subset \mathbb{R}$ and $\Delta$ an admissible bin size,
\begin{eqnarray}\label{smalldiscussion}
&& \sum_{n \in \mathbb{N}_+} {\cal K}(n) \Prob\bigg((x_n,\Delta_n) \in (E \times \{\Delta\}) \,\Big|\,   x_0,\Delta_0 \bigg) \nonumber \\
&& \quad \quad \geq K_{\Delta_0,\Delta} \psi(E,\Delta) \nonumber
\end{eqnarray}
Here $K_{\Delta_0,\Delta}$, denoting a lower bound on the probability of visiting $\Delta$ from $\Delta_0$ in some finite time, is non-zero by (\ref{MoveBin2Bin}) and $\psi$ is positive as the following argument shows: Let $t>0$ be the time stage for which $\Delta_t=\Delta$ and thus by the construction in (\ref{MoveBin2Bin}), with $|h_{t-1}| \leq 1$: $|ax_{t-1}+bu_{t-1}| \leq |a|\Delta_{t-1}/2 = (|a|/\alpha) { \Delta \over 2}$. Thus, it follows that, for $A_1, B_1 \in \mathbb{R}$, $A_1 < B_1$,
\begin{eqnarray}
&& \Prob\bigg(x_{t} \in [A_1, B_1] \,\Big|\,   |ax_{t-1}+bu_{t-1}| \leq |a|\Delta_{t-1}/2, \Delta_{t-1} \bigg) \nonumber \\
&&= \Prob\bigg(ax_{t-1}+bu_{t-1} + d_{t-1} \in [A_1, B_1] \nonumber \\
&& \quad \quad \quad \quad \quad \Big| |ax_{t-1}+bu_{t-1}| \leq |a|\Delta_{t-1}/2, \Delta_{t-1} \bigg) \nonumber \\
&& \geq \min_{|z| \leq { \Delta \over 2} (|a|/\alpha)} \bigg( \Prob(d_{t-1} \in [A_1 - z, B_1 - z] \bigg) > 0
\end{eqnarray}

Now, define the finite set $ C'_{\Delta}:=\{\Delta: L' \leq |\Delta| \leq F', {\log_2(\Delta) \over s} \in \mathbb{N}\}$.
  The chain satisfies the recurrence property that
$\Prob_{(x,\Delta)}(\tau_{C_x \times C'_{\Delta}} < \infty  ) =1$ for any admissible $(x, \Delta)$. This follows, as before in \Section{ProofStocStabSystemTimeDistribution}, from the construction of
\[\Theta_k(\Delta,x) := \Prob(\stp_{1} \geq k\mid x,\Delta), \]
where
\[\stp_{1} = \inf(k > 0: |x_k| \leq 2^{R'-1}\Delta_k, x_0=x, \Delta_0=\Delta )\]
and observing that $\Theta_k(\Delta,x)$ is majorized by a geometric measure with similar steps as in \Section{ProofStocStabSystemTimeDistribution}. Once a state which is perfectly zoomed, that is which satisfies $|x_t| \leq 2^{R'-1}\Delta_t$, is visited, the stopping time analysis can be used to verify that from any initial condition the recurrent set is visited in finite time with probability 1.

In view of (\ref{MoveBin2Bin}), we have that the chain is irreducible.

 We can now show that the set $C_x \times C'_{\Delta}$ is small. We will show first that this set is {\em petite}: A set $D \in\bx$ is petite if there is a probability measure  ${\cal K}$ on the non-negative integers $\mathbb{N}$, and a positive measure $\mu_p$ satisfying $\mu(\state)>0$ and
$$
\sum_{n=0}^{\infty} P^n(x,E) {\cal K}(n) \geq \mu_p(E), \hbox{\it for all \ }  x \in D, \; \mbox{and} \, E
\in \bx .
$$
By Theorem 5.5.7 of \cite{MCSS}, under aperiodicity and irreducibility, every petite set is small. To this end, we will establish aperiodicity at the end of the proof.

To establish the petite set property, we will follow an approach taken by Tweedie \cite{Tweedie2001} which considers the following test, which only depends on the one-stage transition kernel of a Markov chain: If a set $S$ is such that, the following {\em uniform countable additivity} condition
\begin{eqnarray}\label{uniformcountableadditive}
\lim_{n \to \infty} \sup_{x \in S} P(x,B_n) = 0,
\end{eqnarray}
 is satisfied for every sequence $B_n \downarrow \emptyset$, and if the Markov chain is irreducible, then $S$ is petite (see Lemma 4 of Tweedie \cite{Tweedie2001} and Proposition 5.5.5(iii) of Meyn-Tweedie \cite{MCSS}).


Now, the set $C_x \times C'_{\Delta}$ satisfies (\ref{uniformcountableadditive}), since for any given bin size $\Delta'$ in the countable space constructed above, we have that
\begin{eqnarray}
&& \lim_{n \to \infty} \sup_{(x, \Delta) \in C_x \times C'_{\Delta}} \Prob\bigg((x_{t+1},\Delta_{t+1}) \in (B_n \times \Delta') \nonumber \\
&& \quad \quad \quad \quad \quad \quad \quad \Big| x_t=x,\Delta_t=\Delta\bigg) \nonumber \\
&& = \lim_{n \to \infty} \sup_{(x, \Delta) \in C_x \times C'_{\Delta}} \Prob\bigg((ax + bu_t + d_t,\Delta_{t+1}) \nonumber \\
&& \quad \quad \quad \quad \quad \quad \quad \quad \quad \in (B_n \times \Delta') \Big| x_t=x,\Delta_t=\Delta\bigg) \nonumber \\
&& =  \lim_{n \to \infty} \sup_{(x, \Delta) \in C_x \times C'_{\Delta}}\Prob\bigg((d_t,\Delta_{t+1}) \nonumber \\
&& \quad \quad \quad \in \bigg((B_n - (ax + bu_t)) \times \Delta'\bigg) \Big| x_t=x,\Delta_t=\Delta\bigg)\nonumber \\
&& = 0. \nonumber
 \end{eqnarray}
This follows from the fact that the Gaussian random variable $d_1$ satisfies \[\lim_{n \to \infty} \sup_{d_1 \in C_0} \Prob(d_1 \in A_n) = 0,\] uniformly over a compact set $C_0$, for any sequence $A_n \downarrow \emptyset$, since a Gaussian measure admits a uniformly bounded density function.

Therefore, $C_x \times C'_{\Delta}$ is petite.


If the integers $\tilde{A}, \tilde{B}$ are relatively prime, then by B\'ezout's Lemma (see \cite{NumberTheory}), the communication class will include the bin sizes whose logarithms are integer multiples of a constant except those leading to $\Delta < L'$.


We finally show that the Markov chain is aperiodic. This follows from the fact that the smallest admissible state for the quantizer, $\Delta^* = L'$, can be visited in subsequent time stages with non-zero probability, since $$\bigg(\min_{|x| \leq \Delta^*/2} P(d_t \in [-2^{R'-1}\Delta^*-ax, 2^{R'-1}\Delta^*-ax]) \bigg) p >0.$$

\subsubsection{Proof of  Theorem~\ref{Inv}}\label{ProofStocStabSystemProofInv}


With the Lyapunov function $V_0(x_t,\Delta_t) = \log(\Delta^2_t) + B_0$, for $\Delta_{\stp_z }> L$, we have that
\begin{eqnarray}
 && \Expect[V_0(x_{\stp_{z+1}},\Delta_{\stp_{z+1}}) \mid x_{\stp_z},\Delta_{\stp_z}] \nonumber \\
 && = B_0 + \Prob(\stp_{z+1} - \stp_z = 1 \mid x_{\stp_z},\Delta_{\stp_z}) \nonumber \\
&& \quad \quad \quad \quad \times \bigg(2 \log(\alpha) + \log(\Delta_{\stp_z }^2) \bigg) \nonumber \\
&&\quad+ \sum_{k=2}^{\infty} \log(\Delta_{\stp_{z}+k}^2) \Prob(\stp_{z+1} - \stp_z = k\mid x_{\stp_z},\Delta_{\stp_z}) \nonumber
\end{eqnarray}
Thus, the drift satisfies:
\begin{eqnarray}
&& \Expect[V_0(x_{\stp_{z+1}},\Delta_{\stp_{z+1}}) \mid x_{\stp_z},\Delta_{\stp_z}]  -  V_0(x_{\stp_z},\Delta_{\stp_z}) \nonumber \\
&& =  \sum_{k=1}^{\infty} 2 \log( (|a|+\delta)^{(k-1)} \alpha) \Prob(\stp_{z+1} - \stp_z =k \mid x_{\stp_z},\Delta_{\stp_z}) \nonumber \\
&& = 2 \sum_{k=1}^{\infty}  (k-1) \log(|a|+\delta) \Prob(\stp_{z+1} - \stp_z =k \mid x_{\stp_z},\Delta_{\stp_z}) \nonumber \\
&& \quad \quad + 2 \log(\alpha)
. \label{justifyb1}
\end{eqnarray}
By (\ref{definitionofG}),
the summability of $\sum_{k=1}^{\infty} G_k(\Delta_{\stp_z})$,
and the dominated convergence theorem,
\begin{eqnarray}
 && \lim_{\Delta_{\stp_z } \to \infty} \sum_{k=1}^{\infty} (k-1)((1-p)^{k-1} + G_k(\Delta_{\stp_z}) - (1-p)^k) \nonumber \\
 && = \sum_{k=1}^{\infty} \lim_{\Delta_{\stp_z } \to \infty} (k-1) ((1-p)^{k-1} + G_k(\Delta_{\stp_z}) - (1-p)^k) \nonumber \\
 &&  = \sum_{k=1}^{\infty}  p (1-p)^{k-1} (k-1) = p^{-1} - 1
\end{eqnarray}
Provided (\ref{e:Rbdd3}) holds, it follows from Lemma~\ref{KeyBound} and Lemma~\ref{KeyBound2} that for some $b_0>0$,
\begin{eqnarray}\label{e:Rbdd3conc}
&&\lim_{\Delta_{\stp_z } \to \infty} \big\{
  \Expect[V_0(x_{\stp_{z+1}},\Delta_{\stp_{z+1}}) \mid x_{\stp_z},\Delta_{\stp_z}]  -  V_0(x_{\stp_z},\Delta_{\stp_z}) \big\} \nonumber \\
&&=2 \log(\alpha) + 2 \lim_{\Delta_{\stp_z } \to \infty} \big\{ \sum_{k=1}^{\infty} (k-1) \log(|a|+\delta) \nonumber \\ &&  \quad \quad \quad \quad \quad \quad \times \Prob(\stp_{z+1} - \stp_z =k \mid x_{\stp_z},\Delta_{\stp_z})   \big\} \nonumber \\
&& \le -b_0\, .
\end{eqnarray}



For $\Delta_{\stp_z }$ in a compact set and lower bounded by $L'$ defined by (\ref{definitionLprime}), $\Expect[\log(\Delta_{\stp_{z+1}}^2)\mid x_{\stp_z},\Delta_{\stp_z}]$ is uniformly bounded. This follows from the representation of the drift given in (\ref{justifyb1}).
Finally, since,
\[
	G_k(\Delta_{\stp_0}) \leq (1-p)^{k-1} \Gamma_m \Delta_{\stp_0}^{-2m} {1 \over 1 - ((1-p)\xi^{2m})^{-1} },
\]
it follows that $\sum_{k=1}^{\infty} G_k(\Delta_{\stp_0}) k < \infty$ and as a result
\begin{eqnarray}\label{timeBound1}
\sup_{x_{\stp_z},\Delta_{\stp_z}} \Expect_{x_{\stp_z},\Delta_{\stp_z}}[\stp_{z+1} - \stp_{z}] < \infty.
\end{eqnarray}

Consequently, under the bound (\ref{e:Rbdd3}),  there exist $b_0 >0$, $b_1 < \infty, F' >0$ such that,
\begin{eqnarray}
\label{geoDrift}
&& \Expect[V_0(x_{\stp_{z+1}},\Delta_{\stp_{z+1}}) | x_{\stp_z}, \Delta_{\stp_z}] \nonumber \\
&& \quad \quad \leq V_0(x_{\stp_z},\Delta_{\stp_z}) - b_0 + b_1 1_{\{|\Delta_{\stp_z}| \leq F'\}}.
\end{eqnarray}
This combined  with Proposition~\ref{compactSetSmall}, eqn.~(\ref{timeBound1}), and (\ref{geoDrift}), Corollary \ref{corol} leads to positive Harris recurrence.
\qed


\subsubsection{Proof of Theorem~\ref{FiniteMoment}} \label{ProofStocStabSystemFiniteMoment}

First, let us note that by (\ref{MultiBound2}) and (\ref{conditionOnMomentBoundedness}), for every $\kappa >0$, we can find $\Delta_0$ sufficiently large such that
$$
			\lim_{t \to \infty} {\Prob(\stp_1 \geq t | x_0, \Delta_0) \over (1-p+\kappa)^{t-1}} = 0.
$$
Under the assumed bound $(1-p)|a|^2 < 1$,  we can fix $\kappa >0$ such that $(1-p+\kappa) |a+\delta|^{2} < 1.$

Next, observe that for all initial conditions for which $|h_0| \leq 1$,
\begin{eqnarray}\label{Number2}
&&\lim_{\Delta_0 \to \infty} {\Expect[ V_2(x_{\stp_1}, \Delta_{\stp_1})  \mid x_0, \Delta_0] \over V_2(x_{0}, \Delta_{0})}  \nonumber \\
&=&\lim_{\Delta_0 \to \infty} {\Expect[\Delta_{\stp_1}^2 \mid x_0, \Delta_0] \over \Delta_0^2}  \nonumber \\
&=& \lim_{\Delta_0 \to \infty} {1 \over \Delta_0^2} \sum_{k=1}^{\infty} \Prob(\stp_1=k) \Expect[\Delta_k^2|\stp_1=k, x_0, \Delta_0] \nonumber \\
&=& \lim_{\Delta_0 \to \infty}  \alpha^2  \sum_{k=1}^{\infty} \Prob(\stp_1=k) (|a|+\delta)^{2(k-1)} \nonumber \\
&=& p \alpha^2 {1 \over 1 - (1-p) (|a|+\delta)^2 },
\end{eqnarray}
where the last equality follows from Lemma~\ref{KeyBound} and the dominated convergence theorem.

Now, if (\ref{secondMomentBound0}) holds, we can find $\alpha$ such that $R' > \log_2( |a| / \alpha)$, and
\begin{eqnarray}\label{secondMomentBound1}
{p \alpha^2 \over 1 - (1-p)(|a|+\delta)^2} < 1,
\end{eqnarray}
and simultaneously (\ref{e:Rbdd3}) is satisfied. We note that (\ref{secondMomentBound1}) implies (\ref{e:Rbdd3}) since by Jensen's inequality: \[\log(p \alpha^2 +  (1-p)(|a|+\delta)^2) > p \log(\alpha^2) + (1-p) \log((|a|+\delta)^2),\]
and (\ref{e:Rbdd3}) is equivalent to the term on the right hand side being negative.

To establish the required drift equation, we first establish the following bound for all $z \geq 0$:
\begin{eqnarray}\label{keyLemma1}
\kappa \Expect[\sum_{m={\stp_z}}^{\stp_{z+1}-1} x_m^2 \mid x_0,\Delta_0] \leq \Delta_{\stp_z}^2 2^{2(R'-1)},
\end{eqnarray}
for some $\kappa > 0$.

Without loss of generality take $z=0$ so that  $\stp_z=0$.
Observe that
for any $\chi>0$, by H\"older's inequality,
\begin{eqnarray}
&& \Expect[\sum_{t=0}^{\stp_1-1} x_t^2 | x_0,\Delta_0] = \Expect[\sum_{t=0}^{\infty} 1_{\{t < \stp_1\}} x_t^2 | x_0,\Delta_0] \nonumber \\
&& \leq \sum_{t=0}^{\infty} \bigg( \Expect[(1_{\{t < \stp_1\}})^{1 + \chi}| x_0,\Delta_0] \bigg)^{1 \over 1+ \chi} \nonumber \\
&& \quad \quad \quad \quad \quad \quad \times \bigg(\Expect[ x_t^{2({1 + \chi \over \chi})} | x_0,\Delta_0] \bigg)^{\chi \over 1+\chi},
 \label{sonDenklem}
\end{eqnarray}

Moreover,  for some $B_2<\infty$,
\begin{eqnarray} \label{BoundHolder}
&&\Expect[ x_t^{2({1 + \chi \over \chi})} | x_0,\Delta_0] \nonumber \\
&& = \Expect[ a^{2t({1 + \chi \over \chi})}(x_0 + \sum_{i=0}^{t-1} a^{-i-1}d_i)^{2({1 + \chi \over \chi})} | x_0,\Delta_0]  \nonumber \\
&&  = |a|^{2t({1 + \chi \over \chi})} \Expect[  (x_0 + \sum_{i=0}^{t-1} a^{-i-1} d_i))^{2{1 + \chi \over \chi}} | x_0,\Delta_0] \nonumber \\
&& \leq |a|^{2t({1 + \chi \over \chi})} \Expect[  (x_0 + \sum_{i=0}^{\infty} a^{-i-1} d_i)^{2{1 + \chi \over \chi}} | x_0,\Delta_0] \nonumber \\
&& = |a|^{2t({1 + \chi \over \chi})} (2^{R'-1}\Delta_0)^{2{1 + \chi \over \chi}} \nonumber \\
&& \quad \quad \times \Expect[  ({x_0 + \sum_{i=0}^{\infty} a^{-i-1} d_i \over 2^{R'-1}\Delta_0})^{2{1 + \chi \over \chi}} | x_0,\Delta_0] \nonumber \\
&& = |a|^{2t({1 + \chi \over \chi})} (2^{R'-1}\Delta_0)^{2{1 + \chi \over \chi}} \nonumber \\
&& \quad \quad \times \Expect[  (h_0 + {\sum_{i=0}^{\infty} a^{-i-1} d_i \over 2^{R'-1}\Delta_0})^{2{1 + \chi \over \chi}} | x_0,\Delta_0] \nonumber \\
&& < B_2 (2^{R'-1}\Delta_0)^{2{1 + \chi \over \chi}} |a|^{2t({1 + \chi \over \chi})},
  \end{eqnarray}
where the last inequality follows since for every fixed $|h_0| \leq 1$, the random variable $h_0 + (\sum_{i=0}^{\infty} a^{-i-1} d_i) / (2^{R'-1}\Delta_0)$ has a Gaussian distribution with finite moments, uniform on $\Delta_0 \geq L'$. 

Thus,
\begin{eqnarray}
&& \Expect[\sum_{t=0}^{\stp_1-1} x_t^2 | x_0,\Delta_0] \nonumber \\
 && \leq \sum_{t=0}^{\infty} \bigg(\Expect[(1_{\{t < \stp_1\}})^{1 + \chi}| x_0,\Delta_0] \bigg)^{1 \over 1+ \chi} \nonumber \\
&& \quad \quad \quad \quad \quad \quad \times \bigg(B_2 (2^{R'-1}\Delta_0)^{2{1 + \chi \over \chi}}  |a|^{2t({1 + \chi \over \chi})} \bigg)^{\chi \over 1+\chi} \nonumber \\
&&= (2^{R'-1}\Delta_0)^2 \sum_{t=0}^{\infty} \bigg( \Prob(\stp_1 \geq t+1 | x_0, \Delta_0) \bigg)^{1 \over 1+ \chi} \nonumber \\
&& \quad \quad \quad \quad \quad \quad \times \bigg( B_2 |a|^{2t({1 + \chi \over \chi})} \bigg)^{\chi \over 1+\chi} \nonumber \\
&&< \zeta_{B_2} (2^{R'-1}\Delta_0)^2 \nonumber 
\end{eqnarray}
for some $\zeta_{B_2} < \infty$.

The last inequality is due to the fact there exists $\kappa >0$ such that $$\lim_{t \to \infty} {\Prob(\stp_1 \geq t | x_0, \Delta_0) \over (1-p+\kappa)^{t-1}} = 0,$$
and we can pick $\chi >0$ with $(1-p+\kappa) |a|^{2(1 + \chi)} < 1.$ Such $\chi $ and $\kappa$ exist by the hypothesis that $(1-p)|a|^2 < 1$.

Hence, with $0 < \epsilon < 1 - p \alpha^2 / [ 1 - (1-p)(|a|+\delta)^2]$,
\[
  \delta(x,\Delta) = \epsilon \Delta^2, \quad f(x,\Delta) = {\epsilon \over \zeta_{B_2} 2^{2(R'-1)}} x^2,
\]
$C$ a compact set, and $V_2(x,\Delta)=\Delta^2$, \Theorem{thm5} applies and $\lim_{t \to \infty} E[x_t^2] < \infty$.
\qed

\subsubsection{Proof of Theorem~\ref{FinitemMoment}}\label{ProofStocStabSystemFinitemMoment}

The proof follows closely that of Theorem~\ref{FiniteMoment}.


Again by H\"older's inequality, for any
  $\chi>0$,
\begin{eqnarray}
&& \Expect\bigg[\sum_{t=0}^{\stp_1-1} |x_t|^m \,\Big|\,   x_0,\Delta_0\bigg] = \Expect\bigg[\sum_{t=0}^{\infty} 1_{\{t < \stp_1\}} |x_t|^m \,\Big|\,   x_0,\Delta_0\bigg] \nonumber \\
&& \leq \sum_{t=0}^{\infty} \bigg( \Expect[(1_{\{t < \stp_1\}})^{1 + \chi}| x_0,\Delta_0] \bigg)^{1 \over 1+ \chi} \nonumber \\
&& \quad \quad \quad \quad \quad \quad \times \bigg(\Expect[ |x_t|^{m({1+ \chi \over  \chi})} | x_0,\Delta_0] \bigg)^{ \chi \over 1+ \chi}.
 \label{sonDenklem}
\end{eqnarray}
As in (\ref{BoundHolder}),  for some $B_m < \infty$,
\[
\Expect[ |x_t|^{m({1+ \chi \over  \chi})} | x_0,\Delta_0] \leq B_m (\Delta_0^m 2^{m(R'-1)})^{1+ \chi \over  \chi} |a|^{mt({1+ \chi \over  \chi})}
\]
 %
and consequently,
\begin{eqnarray}
&& \Expect\big[\sum_{t=0}^{\tau_1-1} |x_t|^m | x_0,\Delta_0\big] \nonumber \\
&& \leq \Delta_0^m 2^{m(R'-1)}  \nonumber \\
&& \quad \times \sum_{t=0}^{\infty} \bigg( \Prob(\stp_1 \geq t+1 | x_0,\Delta_0) \bigg)^{1 \over 1+ \chi} \bigg( B_m^{\chi \over 1+\chi} |a|^{mt} \bigg) \nonumber \\
&& < \zeta_{B_m} (2^{R'-1}\Delta_0)^m
\end{eqnarray}
where, once again, the last inequality is due to the fact that there exists a $\kappa >0$ such that $$\lim_{t \to \infty} {\Prob(\stp_1 \geq t | x_0, \Delta_0) \over (1-p+\kappa)^{t-1}} = 0,$$
and we can pick $\chi >0$ such that $(1-p+\kappa) |a|^{m(1 + \chi)} < 1$;
such $\chi $ and $\kappa$ exist by the property $(1-p)|a|^m < 1$.



Hence, with $0 < \epsilon < 1 - p \alpha^m /[ 1 - (1-p)(|a|+\delta)^m]$,
\[
  \delta(x,\Delta) = \epsilon \Delta^m, \quad f(x,\Delta) = {\epsilon \over \zeta_{B_m} 2^{2(R'-1)}} |x|^m,
\]
$C$ a compact set,
and $V_m(x,\Delta)=\Delta^m$, \Theorem{thm5} applies, establishing in the desired conclusions, and in particular that $\lim_{t \to \infty} E[|x_t|^m]  $ exists and is finite.
\qed

\section{Acknowledgements}
Financial support from the AFOSR grant FA9550-09-1-0190  is gratefully acknowledged.

We are grateful to three reviewers and the associate editor for comments which have led to significant improvements in the presentation of the paper.
\def\cprime{$'$} 

\begin{biography}
{\bf Prof. Serdar Y\"uksel} (S'02) was born in Lice, Turkey in 1979. He received his BSc degree in Electrical and
Electronics Engineering from Bilkent University in 2001; MS and PhD degrees in Electrical and
Computer Engineering from the University of Illinois at Urbana-Champaign in 2003 and 2006,
respectively. He was a post-doctoral researcher at Yale University before joining
Queen's University at the Department of Mathematics and Statistics in 2007. His research interests are on stochastic control, decentralized control, information theory, source coding theory, multi-terminal control and communication systems and stochastic processes.
\end{biography}

\begin{biography}
{\bf Prof. Sean P. Meyn} received the B.A. degree in mathematics from UCLA in 1982, and the Ph.D. degree in electrical engineering from McGill University in 1987 (with Prof. P. Caines).  After 22 years as a professor at the University of Illinois,  he is now Robert C. Pittman Eminent Scholar Chair in the Dept. of ECE at the University of Florida, and director of the new Laboratory for Cognition \& Control.  His research interests include stochastic processes, optimization, complex networkks, information theory, and power and energy systems.\end{biography}

\end{document}